\documentclass[11pt,russian]{article}
\usepackage{amssymb,amsmath,comment}
\usepackage[cp1251]{inputenc}
\usepackage[russian]{babel}

\def\risk{\mathit{R}}

\title{Экспоненциальное взвешивание и оракульные неравенства для проекционных оценок\thanks{Работа выполнена при поддержке 
Лаборатории структурных методов анализа данных в предсказательном моделировании, МФТИ, грант правительства РФ дог. 11 11.G34.31.0073
}}
\author{Голубев Г.К. \thanks{CNRS,
Институт Проблем Передачи Информации РАН 
и Московский Физико-Технический Институт,\quad 
e-mail: {\tt golubev.yuri@gmail.com}}}

\date{}
\newtheorem{theorem}{Теорема}
\newtheorem{lemma}{Лемма}

\begin{document}
\renewcommand{\figurename}{Рис.}
\renewcommand{\refname}{\centerline{СПИСОК ЛИТЕРАТУРЫ}}
\renewcommand{\abstractname}{}

 \maketitle
\begin{abstract}
Рассматривается задача оценивания неизвестного вектора,
наблюдаемого на фоне белого гауссовского шума. Для  оценки этого
вектора используется семейство проекционных оценок и задача
состоит в том, чтобы на основе наблюдений выбрать наилучшую выпуклую 
комбинацию оценок из этого семейства. В работе изучается метод построения
 оценок, основанный на так называемом экспоненциальном взвешивании и приводится верхняя граница для среднеквадратичного риска этого метода.

\end{abstract}

\section{Введение и основные результаты}

В настоящей работе рассматривается   простейшая  линейная   модель, в
которой нужно оценить неизвестный вектор
 $\mu \in {l}_2(1,\infty)$  на основе  наблюдений
\begin{equation}\label{equ.1}
Y_{i}=\mu _{i}+\sigma \xi _{i} , \quad i=1,2,\ldots,
\end{equation}
где $\xi _{i}$ -- белый гауссовский шум, т. е. $\xi _{i}$ являются
независимыми гауссовскими  случайными величинами с нулевым средним и единичной дисперсией (${\cal N}(0,1)$). Далее, чтобы упростить многочисленные
технические детали, предполагается, что параметр $\sigma>0$
известен.

Для краткости, будем обозначать векторы
$(\mu_1,\mu_2,\ldots)^\top$  и $(Y_1,Y_2,\ldots)^\top$ как
 $\mu$ и $Y$. Пусть $\hat{\mu
}(Y)=(\hat{\mu}_1(Y),\hat{\mu}_2(Y),\ldots)^\top$ --
некоторая оценка вектора $\mu$. Ее риск будем измерять
следующей величиной
\[
\risk(\hat{\mu },\mu )=\mathbf{E}_{\mu }\Vert \hat{\mu }(Y)
-\mu \Vert ^{2},
\]
здесь и далее $\mathbf{E}_{\mu }$ -- математическое ожидание по
мере $\mathbf{P}_{\mu}$, порожденной наблюдениями (\ref{equ.1}), а
$\left\| \cdot \right\| $ и $\langle \cdot ,\cdot \rangle$
обозначают норму и скалярное произведение в $l_2(1,\infty)$
\[
\left\| x\right\| ^{2}=\sum_{i=1}^{\infty }x_{i}^{2},\quad \langle
x,y\rangle =\sum_{i=1}^{\infty }x_{i}y_{i}.
\]

Для оценивания неизвестного вектора $\mu$ будем использовать
проекционные оценки
\[
\hat\mu _{i}^{m}=\mathbf{1}\{i\le m\}Y_{i},
\quad 
m\in \mathcal{M};
\]
здесь $\mathcal{M}$ -- некоторое  ограниченное множество целых чисел. Как правило, в качестве $\mathcal{M}$  используется $\{1,2,\ldots,n\}$.  

Имея в своем распоряжении это семейство оценок, мы будем оценивать
    $\mu$ с помощью выпуклой комбинации оценок $\hat\mu^m,\ m\in \mathcal{M}$
\begin{equation*}
\bar{\mu}^w(Y)=\sum_{m\in \mathcal{M}} w_m(Y)\hat\mu^m(Y),
\end{equation*}
где веса $w_m(Y)$ зависят от наблюдений, положительны и таковы, что
\[
\sum_{m\in \mathcal{M}} w_m(Y)=1.
\]

Задача, рассматриваемая в этой статье, состоит в том, чтобы найти веса $w_m(Y)$, которые минимизируют риск оценки  $\hat{\mu}^w(Y)$.

По-видимому, первые подходы к решению этой задачи  связаны с  атомарными весами и  с идеей несмещенного оценивания
рисков оценок $\hat{\mu}^m$  \cite{A}.
Современная математическая литература в этой области настолько
обширна, что ее даже не имеет  смысла здесь цитировать, но нельзя не упомянуть статью \cite{K},
которая является одной из классических работ по непараметрическому оцениванию.
В подходе, основанном на несмещенном оценивании риска,  выбирается проекционный метод, имеющий минимальную несмещенную оценку риска. Точнее,  определим
\[
\hat{m}(Y)=\arg \min_{m\in \mathcal{M}} \bigl\{\bar{r}(Y,\mu^m)
\bigr\},
\]
где
\begin{equation}\label{equ.2}
\bar{r}(Y,\hat{\mu}^m)=-\sum_{i=1}^m Y_i^2 +2\sigma^2 m
\end{equation}
-- несмещенные оценки рисков оценок $\hat{\mu}^m$. Заметим, что, на самом деле, величины $\bar{r}(Y,\hat{\mu}^m)$ являются несмещенными оценками риска с точностью до аддитивной постоянной. Точнее,  несмещенная оценка риска компоненты $\hat{\mu}^m_i$ имеет  следующий вид:
\[
(Y_i-\hat{\mu}^m_i)^2+2\sigma^2\mathbf{1}\{i\le m\}-\sigma^2. 
\] 

Для метода, основанного на несмещенном оценивании риска,  справедлив следующий факт, который нетрудно вывести из \cite{K}.
\begin{theorem}\label{th.1} Пусть  $w_m(Y)=\delta(m,\hat{m}(Y))$, где $\delta(x,x)=1$ и $\delta(x,y)=0,\ x\neq y$. Тогда
для риска оценки $\bar{\mu}^w(Y)$ справедлива следующая верхняя граница:
\begin{equation*}
R(\bar\mu^w,\mu)\le r^{\mathcal{M}}(\mu)+K\sigma^2\sqrt{\frac{r^{\mathcal{M}}(\mu)}{\sigma^2}},
\end{equation*}
где
\begin{equation}\label{equ.3}
r^{\mathcal{M}}(\mu)=\min_{m\in \mathcal{M}} \biggl\{\sum_{i=m+1}^\infty\mu_i^2+\sigma^2 m \biggr\},
\end{equation}
и $K$ -- некоторая универсальная постоянная.
\end{theorem}

Величина $r^{\mathcal{M}}(\mu)$ часто называется риском оракула. Действительно, если мы предположим, что у нас  имеется доступ к оракулу, который для любой оценки 
$\hat{\mu}$ может точно предсказывать ее риск $\mathbf{E}_\mu \|\mu-\hat{\mu}\|^2$, то $r^{\mathcal{M}}(\mu)$ это минимальный риск, который может достигнут при помощи оракула. В действительности, мы, конечно,  проиграем оракулу и формула (\ref{equ.3}) показывает, что величина этого проигрыша (в случае проекционных
оценок)  не больше, чем $K\sigma^2\sqrt{r^{\mathcal{M}}(\mu)/\sigma^2}$. Чтобы понять насколько это хорошо или плохо, рассмотрим два случая:
\begin{enumerate}
\item Малая эффективная размерность вектора $\mu$, т.е.  $r^{\mathcal{M}}(\mu)\approx \sigma^2$. 
\item Большая эффективная размерность вектора $\mu$, т.е.  $r^{\mathcal{M}}(\mu)\gg \sigma^2$. 
\end{enumerate}    
 В первом случае риск нашего метода будет  иметь порядок $(K+1)r^{\mathcal{M}}(\mu)$, то
есть он превосходит в $K+1$ раз риск оракула. В случае же когда эффективная размерность  $\mu$ велика,
отношение риска нашего метода и риска оракула приближается к 1. Поэтому часто эту ситуацию называют адаптивным оцениванием.            


Другой хорошо известный результат \cite{LB}  связан с так называемым экспоненциальным взвешиванием. Положим
\begin{equation}\label{equ.4}
w^*_m(Y)=\exp\biggl(-\frac{\bar{r}(Y,\hat{\mu}^m)}{4\sigma^2}\biggr)\biggl/
\sum_{s\in \mathcal{M}} \exp\biggl(-\frac{\bar{r}(Y,\hat{\mu}^s)}{4\sigma^2}\biggr).
\end{equation}
\begin{theorem}\label{th.2} Для риска оценки $\bar\mu^{w^*}$ справедлива следующая верхняя граница:
\begin{equation}\label{eq.5}
R(\bar\mu^{w^*},\mu)\le r^{\mathcal{M}}(\mu)+4\sigma^2\log(\#\mathcal{M}),
\end{equation}
где
$\#\mathcal{M}$ обозначает число элементов множества $\mathcal{M}$.  
\end{theorem}

Идея использовать неатомарные  веса имеет довольно долгую историю. По-видимому,
первые математически строгие результаты в этой области были получены А.С. Немировским (см., например, \cite{N}) в предположении, что имеется дополнительная обучающая выборка, которая используется для выбора весов. 
Позднее, результаты близкие к теореме \ref{th.2} были получены в серии статей О. Катони (см., например, монографию \cite{C}, которая суммирует эти результаты).
Из недавних работ читателям, интересующимся этим методом,  можно обратить внимание на \cite{DS} и \cite{RT}. По существу дела, эти работы  близки к \cite{LB} и обобщают результаты этой статьи на  более сложные статистические модели. 

Первый вопрос, который возникает при взгляде на теоремы \ref{th.1} и \ref{th.2} это, естественно, вопрос о том помогают ли эти результаты понять, какой же метод
лучше. К сожалению, этого сделать фактически нельзя. Дело в том, что верхняя граница в теореме \ref{th.2} стремиться к бесконечности при $\#\mathcal{M}\rightarrow\infty$ и  в зависимости от величины $\#\mathcal{M}$ и вектора $\mu$  граница (\ref{eq.5}) может быть как лучше, так и хуже
(\ref{equ.3}).  Понять из теорем  \ref{th.1} и \ref{th.2} причину этого эффекта невозможно: это может происходить как от того, что  граница
(\ref{eq.5}) неоправданно завышена, так и от того, что  методы могут оказаться принципиально несравнимыми.

Цель настоящей работы улучшить верхнюю границу (\ref{eq.5}) так, чтобы она оставалась   ограниченной при неограниченном увеличении  числа элементов множества $\mathcal{M}$. Точнее, мы покажем, что справедлив следующий результат.

\begin{theorem}\label{th.3} 
Для риска оценки $\bar\mu^{w^*}$ справедлива следующая верхняя граница:
\begin{equation}\label{equ.6}
R(\bar\mu^{w^*},\mu)\le r^{\mathcal{M}}(\mu)+4\sigma^2\log\biggl\{\frac{r^{\mathcal{M}}(\mu)}{\sigma^2}\biggl[1+\Psi\biggl(\frac{\sigma^2}{r^{\mathcal{M}}(\mu)}\biggr)\biggr]\biggr\},
\end{equation}
где $\Psi(r),\ r\in [0,1]$ -- некоторая ограниченная функция такая, что $$\lim_{r\rightarrow 0}\Psi(r)=0.$$
\end{theorem}

Доказательство этого результата проводится  с помощью  комбинации методов из \cite{LB} и \cite{G} и основано на том, что начиная с некоторого (случайного) $m^\circ$, несмещенные оценки рисков $\bar{r}(Y,\hat{\mu}^m)$ ограничены снизу линейной функцией от $m$.
Статистический смысл теоремы \ref{th.3} достаточно прозрачен. Неравенство (\ref{equ.6}) говорит том, что экспоненциальное взвешивание уменьшает плату за отсутствие оракула. Если в классическом методе несмещенного оценивания риска мы должны платить $K\sigma^2\sqrt{r^{\mathcal{M}}(\mu)/\sigma^2}$, то при
экспоненциальном взвешивании плата ограничена, грубо говоря, величиной $4\sigma^2 \log[r^{\mathcal{M}}(\mu)/\sigma^2]$ при больших отношениях $r^{\mathcal{M}}(\mu)/\sigma^2\gg 1$.  При малых отношениях $r^{\mathcal{M}}(\mu)/\sigma^2$ границы из теорем \ref{th.1} и \ref{th.3} становятся эквивалентными, т.к. $r^{\mathcal{M}}(\mu)/\sigma^2\ge 1$.  Это подтверждает гипотезу о том, что экспоненциальное взвешивание является более эффективным, чем классические методы выбора моделей, основанные на несмещенном оценивании риска. 

Заметим также, что в качестве границы для риска $R(\bar\mu^{w^*},\mu)$ можно, конечно, брать минимум из правых частей (\ref{eq.5}) и (\ref{equ.6}). 
\section{Доказательства}
Отправной точкой для доказательства теоремы \ref{th.3} является следующий факт:
\begin{lemma} \label{lemma.1}Пусть $\|\mu\|<\infty$.  
Тогда
\begin{equation*}
R(\bar\mu^{w^*},\mu)\le r^{\mathcal{M}}(\mu) +\mathbf{E}\sum_{m\in\mathcal{M}} w^*_m(Y)[\bar{r}(Y,\hat{\mu}^m)-\bar{r}^{\mathcal{M}}(Y)];
\end{equation*}
где ${r}^{\mathcal{M}}(\mu)$ определено в (\ref{equ.3}), а $\bar{r}^{\mathcal{M}}(Y)=\min_{m\in \mathcal{M}}\bar{r}(Y,\hat{\mu}^m)$.
\end{lemma}
\textit{Доказательство.} Оно приводится для удобства читателей и с точностью до второстепенных деталей следует \cite{LB}. Начнем с хорошо известной формулы Стейна \cite{St} для несмещенного оценивания риска.
Предположим, что для оценивания вектора $\mu$ используется нелинейная оценка $\hat{\mu}(Y)$, компоненты
которой и имеют вид
\begin{equation*}
\hat{\mu}_i(Y)=Y_i+\phi_i(Y),
\end{equation*}   
где $\phi_i(\cdot)$ -- некоторая  дифференцируемая функция.  
Тогда, интегрируя по частям, находим
\begin{equation*}
\begin{split}
\mathbf{E}(\mu_i-\hat{\mu}_i)^2=\mathbf{E}[\sigma \xi_i+\phi_i(\mu+\sigma\xi_i)]^2=\sigma^2+2\sigma\mathbf{E}\xi_i \phi_i(\mu+\sigma\xi)+ \mathbf{E}\phi_i^2(Y)\\
 =\sigma^2 +2\sigma^2 \mathbf{E}\frac{\partial \phi_i(Y)}{\partial Y_i} + \mathbf{E}\phi_i^2(Y).
\end{split}
\end{equation*}
Другими словами
\begin{equation}\label{equ.7}
\mathbf{E}(\mu_i-\hat{\mu}_i)^2= \mathbf{E}\bar{R}(Y_i,\hat{\mu}_i),
\end{equation} 
где 
\[
\bar{R}(Y_i,\hat{\mu}_i)=(Y_i-\hat{\mu}_i)^2+2\sigma^2 \frac{\partial \hat{\mu}_i(Y)}{\partial Y_i} +\sigma^2.
\]
Поэтому величину $\bar{R}(Y_i,\hat{\mu}_i)$ называют несмещенной оценкой риска оценки $\hat{\mu}_i$.

Применим теперь эти формулы к семейству  проекционных оценок $\hat{\mu}^m(Y),\ m\in \mathcal{M}$ и их выпуклой комбинации
\[
\bar{\mu}^{w^*} =\sum_{m \in \mathcal{M}} w_m^*(Y)\hat{\mu}^m(Y)
\]   
с весами из (\ref{equ.4}).
 Пусть как и ранее
\[
\bar{r}(Y,\hat{\mu}^m) = \|\hat{\mu}^m\|^2-2\langle\hat{\mu}^m,Y\rangle +2\sigma^2\sum_{i=1}^{\max(\mathcal{M})} \frac{\partial\hat{\mu}^m_i(Y)}{\partial Y_i}
\] 
несмещенные оценки для рисков  оценок $\hat{\mu}^m$ с точностью до аддитивной постоянной.
Здесь  и далее $\max(\mathcal{M})$ обозначает максимальный элемент в $\mathcal{M}$. 

Заметим, что с учетом этого определения  все суммы, входящие в  несмещенные 
оценки рисков $\hat{\mu}^m(Y)$ и  $\hat{\mu}^{w^*}(Y)$,  становятся конечными поскольку множество $\mathcal{M}$ предполагается ограниченным.
  
 Наша цель связать среднеквадратичный риск оценки 
 $\bar{\mu}^{w^*}$ и  величины $\bar{r}(Y,\hat{\mu}^m)$.
Заметим, что
\begin{equation}\label{equ.8}
\begin{split}
\|\bar{\mu}^{w^*}\|^2-2\langle \bar{\mu}^{w^*} Y\rangle= \\ =\sum_{m \in \mathcal{M}} w_m^*(Y)\Bigl\{\|\hat{\mu}^m-\hat{\mu}^m+\bar{\mu}^{w^*}\|^2-2
\langle\hat{\mu}^m-\hat{\mu}^m+\bar{\mu}^{w^*},Y\rangle\Bigr\}\\
=\sum_{m \in \mathcal{M}} w_m^*(Y)\Bigl\{\|\bar{\mu}^{w^*}-\hat{\mu}^m\|^2+\|\hat{\mu}^m\|^2-2\langle \hat{\mu}^m ,Y\rangle^2\\
+2\langle Y-\hat{\mu}^m,\hat{\mu}^m-\bar{\mu}^{w^*}\rangle
\Bigr\}\\
=\sum_{m \in \mathcal{M}} w_m^*(Y)\Bigl\{\|\bar{\mu}^{w^*}-\hat{\mu}^m\|^2+\|
\hat{\mu}^m\|^2-2\langle \hat{\mu}^m, Y\rangle \\ +2\langle \bar{\mu}^{w^*}-\hat{\mu}^m,\hat{\mu}^m-\bar{\mu}^{w^*}\rangle
\Bigr\}\\
=\sum_{m \in \mathcal{M}} w_m^*(Y)\Bigl\{\|\hat{\mu}^m\|^2-2\langle \hat{\mu}^m ,Y\rangle \Bigr\}-\sum_{m\in\mathcal{M} } w_m^*(Y)\|\bar{\mu}^{w^*}-\hat{\mu}^m\|^2.
\end{split}
\end{equation}

Из определения оценки $\bar{\mu}^{w^*}$ и  из (\ref{equ.7}) мы получаем, что
\begin{equation*}
\begin{split}
\mathbf{E}\|\bar{\mu}^{w^*}-\mu\|^2  
= \mathbf{E}\sum_{i=1}^{\max(\mathcal{M})}\biggl[(\hat{\mu}_i^{w^*})^2-2\hat{\mu}_i^{w^*}Y_i+2\sigma^2\sum_{m \in \mathcal{M}} w_m^*(Y)\frac{\partial \hat{\mu}_i^m(Y)}{\partial Y_i}\\
+2\sigma^2\sum_{m \in \mathcal{M}} \hat{\mu}_i^m(Y)\frac{\partial w_m^*(Y)}{\partial Y_i} \biggr]+\|\mu\|^2.
\end{split}
\end{equation*}
Воспользовавшись  (\ref{equ.8}), продолжим это соотношение следующим образом: 
\begin{equation}\label{equ.9}
\begin{split}
\mathbf{E}\|\bar{\mu}^{w^*}-\mu\|^2=\|\mu\|^2+\mathbf{E}
\sum_{m \in \mathcal{M}} w_m^*(Y)\Bigl\{\|\hat{\mu}^m\|^2-2\langle \hat{\mu}^m,Y\rangle\Bigr\}\\ -\sum_{m \in \mathcal{M}} w_m^*(Y)\|\bar{\mu}^{w^*}-\hat{\mu}^m\|^2
+2\sigma^2\sum_{m \in \mathcal{M}} w_m^*(Y)\frac{\partial \hat{\mu}^m(Y)}{\partial Y_i}\\
=\|\mu\|^2+\mathbf{E}
\sum_{m \in \mathcal{M}} w_m^*(Y)\Bigl\{\|\hat{\mu}^m\|^2-2\langle \hat{\mu}^m,Y\rangle\Bigr\}\\ -\sum_{m \in \mathcal{M}} w_m^*(Y)\|\bar{\mu}^{w^*}-\hat{\mu}^m\|^2\\
+2\sigma^2\sum_{m \in \mathcal{M}}\sum_{i=1}^{\max(\mathcal{M})} w_m^*(Y)\frac{\partial \hat{\mu}_i^m(Y)}{\partial Y_i}
+2\sigma^2\sum_{m \in \mathcal{M}} \sum_{i=1}^{\max(\mathcal{M})}\hat{\mu}_i^m(Y)\frac{\partial w_m^*(Y)}{\partial Y_i} \biggr]\\
=\|\mu\|^2+ \mathbf{E}
\sum_{m \in \mathcal{M}} w_m^*(Y)\biggl\{ \bar{r}(Y,\hat{\mu}^m)-\|\bar{\mu}^{w^*}-\hat{\mu}^m\|^2
\\ +2\sigma^2 \sum_{i=1}^{\max(\mathcal{M})}\hat{\mu}_i^m(Y)\frac{\partial\log[ w_m^*(Y)]}{\partial Y_i} \biggr\}.
\end{split}
\end{equation} 
Далее, так как $\sum_{m \in \mathcal{M}} w_m^*(Y)=1$,  то
\[
\sum_{m \in \mathcal{M}} w_m^*(Y) \frac{\partial\log[ w_k^*(Y)]}{\partial Y_i}=0
\]
и из (\ref{equ.9}) имеем
\begin{equation}\label{equ.10}
\begin{split}
\mathbf{E}\|\bar{\mu}^{w^*}-\mu\|^2=
\mathbf{E}
\sum_{m \in \mathcal{M}} w_m^*(Y)\biggl\{ \bar{r}(Y,\hat{\mu}^m)-\|\bar{\mu}^{w^*}-\hat{\mu}^m\|^2\\
+2\sigma^2 \sum_{i=1}^{\max(\mathcal{M})} (\hat{\mu}_i^m(Y)-\hat{\mu}_i^{w^*})\frac{\partial\log[ w_m^*(Y)]}{\partial Y_i} \biggr\}.
\end{split}
\end{equation}

Наш следующий  шаг -- вычислить частные производные в правой части этого  равенства. Заметив, что
\begin{equation*}
\begin{split}
\frac{\partial \bar{r}[Y,\hat{\mu}^m] }{\partial Y_i}=
\frac{\partial \phantom{Y_i} }{\partial Y_i}\bigl[\|\hat{\mu}^m\|^2-2\langle \hat{\mu}^m  ,Y\rangle  +2\sigma^2 m\bigr]=
-2\hat{\mu}^m_i,
\end{split}
\end{equation*}
находим 
\begin{equation*}
\begin{split}
\frac{\partial\log[ w_k^*(Y)]}{\partial Y_i} =-
\frac{1}{4\sigma^2}\frac{\partial \bar{r}[Y,\hat{\mu}^k] }{\partial Y_i}
+
\frac{1}{4\sigma^2}\sum_{s \in \mathcal{M}} w_s^*(Y)\frac{\partial \bar{r}[Y,\hat{\mu}^s] }{\partial Y_i}\\
=\frac{\hat{\mu}^k_i}{2\sigma^2} -\frac{1}{2\sigma^2}\sum_{s \in \mathcal{M}}
w_s^*(Y)\hat{\mu}^s_i=\frac{\hat{\mu}^k_i-\hat{\mu}_i^{w^*}}{2\sigma^2}.
\end{split}
\end{equation*}
Подставляя это соотношение в (\ref{equ.10}), мы получаем
\begin{equation*}
\begin{split}
\mathbf{E}\|\bar{\mu}^{w^*}-\mu\|^2=\|\mu\|^2+\mathbf{E}  \bar{r}^{\mathcal{M}}(Y)+
\mathbf{E}
\sum_{m \in \mathcal{M}} w_k^*(Y)\bigl[ \bar{r}(Y,\hat{\mu}^m)- \bar{r}^{\mathcal{M}}(Y)\bigr].
\end{split}
\end{equation*}
Чтобы завершить доказательство теоремы осталось заметить, что
\[
\mathbf{E}  \bar{r}^{\mathcal{M}}(Y)\le \min_{m\mathcal{M}}\mathbf{E}   \bar{r}(Y,\hat{\mu}^m)=r^{\mathcal{M}}(\mu)-\|\mu\|^2.
\]

\quad $\blacksquare$

\bigskip
 Далее нам потребуются также  простые и хорошо известными вероятностные  факты, которые собраны в следующей лемме.

\begin{lemma}
\label{lemma.2} Пусть $\xi_i$ -- независимые $\mathcal{N}(0,1)$.
Тогда
\begin{eqnarray}\label{equ.11}
\mathbf{E}\max_{k\geq 1}\biggl\{ \sum_{i=1}^{k}(\xi
_{i}^{2}-1)-U(\alpha )k\biggr\} \leq \frac{1}{\alpha },\\ \label{equ.12}
\mathbf{E}\max_{k\geq 1}\left\{ \sum_{i=k}^\infty \xi_i\mu_i-\frac{\alpha }{2}\sum_{i=k}^\infty \mu_i^2\right\} \leq \frac{1%
}{\alpha }, \\
\mathbf{E}\max_{k\geq 1}\biggl\{ \sum_{i=1}^{k}(1-\xi
_{i}^{2})-U^*(\alpha) k\biggr\} \leq \frac{1}{\alpha }, \label{equ.13}
\end{eqnarray}
где 
\[
U(\alpha)=-\frac{\alpha+\log(1-2\alpha)/2}{\alpha}, \quad U^*(\alpha)=\frac{\alpha-\log(1+2\alpha)/2}{\alpha}.
\]
\end{lemma}
\textit{Доказательство.} 
 Обозначим для краткости
\[
\varphi (\alpha )=\mathbf{E}\exp [\alpha (\xi _{1}^{2}-1)]=\frac{\exp
(-\alpha )}{\sqrt{1-2\alpha }}.
\]
Несложно проверить, что случайный процесс
\[
m_{k}=\exp \biggl\{ \alpha \sum_{i=1}^{k}(\xi _{i}^{2}-1)\biggr\}
\varphi ^{-k}(\alpha )
\]
является мартингалом и, следовательно,  $
\mathbf{E}m_{\tau }=1
$ для любого момента остановки $\tau $ такого, что  $\mathbf{E}%
\tau <\infty$.

Рассмотрим следующий момент остановки
\[
\tau (x,A)=\min \{k\ge 0:S_{k}\notin [x,-A]\},
\]
где
\[
S_0=0, \quad S_{k}=\sum_{i=1}^{k}(\xi _{i}^{2}-1)-U(\alpha )k,\ k\ge 1.
\]
Из тождества $\mathbf{E}m_{\tau (x,A)}=1$ находим
\begin{equation*}
\begin{split}
1 &\geq \mathbf{E}m_{\tau (x,A)}\mathbf{1}%
\{S_{\tau (x,A)}\geq x\} \\ &=\mathbf{E}\exp \biggl\{ \alpha\biggr[
\sum_{i=1}^{\tau (x,A)}(\xi _{i}^{2}-1)-\tau (x,A)\alpha^{-1}\log \varphi
(\alpha )\biggr]\biggr\}   \mathbf{1}\{S_{\tau (x,A)}\geq x\}
\\ &=\mathbf{E}\exp \biggl\{ \alpha\biggr[
\sum_{i=1}^{\tau (x,A)}(\xi _{i}^{2}-1)-\tau (x,A)U(\alpha)\biggr]\biggr\}   \mathbf{1}\{S_{\tau (x,A)}\geq x\}
 \\ &\geq \exp (\alpha
x)\mathbf{E1}\{S_{\tau (x,A)}\geq x\}.
\end{split}
\end{equation*}
Следовательно,
\[
\mathbf{P}\{S_{\tau (x,A)}\geq x\}\leq \exp (-\alpha x).
\]
Наконец, используя следующую формулу
\[
\mathbf{P}\{\max_{k\geq 1}S_{k}>x\}=\lim_{A\rightarrow \infty }\mathbf{P}%
\{S_{\tau (x,A)}\geq x\},
\]
получаем, что
\[
\mathbf{E}\max_{k\geq 1}S_{k}=\int_{0}^{\infty
}\mathbf{P}\{\max_{k\geq 1}S_{k}>x\}dx\leq \frac{1}{\alpha
},
\]
таким образом, доказывая (\ref{equ.11}). 
Для доказательства  (\ref{equ.12}) и (\ref{equ.13}) можно применить совершенно аналогичные
рассуждения.  \quad $\blacksquare$

\begin{lemma} \label{lemma.3} Пусть $U^{-1}(\alpha)$ и $U^{*-1}(\alpha)$ -- функции обратные к $U(\alpha)$ и $U^{-1}(\alpha)$ соотвественно. Для них справедливы следующие неравенства:
\begin{align}
\label{equ.14}
U^{-1}(\alpha)\ge \frac{\alpha}{1+2\alpha}, \quad
U^{*-1}(\alpha)\ge \alpha.
\end{align}

\end{lemma}
\textit{Доказательство.} Заметим, что функция
\begin{equation*}
f(x)\stackrel{\rm def}{=}\log(1-x)-x-\frac{x^2}{2(1-x)},\ x\in [0,1),
\end{equation*}
отрицательна. 
Для проверки этого факта достаточно взглянуть на первую производную $f(x)$
$$
f'(x)=-\frac{x^2}{(1-x)^2}.
$$ 
Поэтому $f(x)$ убывает при $x\ge 0$ и достигает максимума в точке $x=0$. При этом $f(0)=0$. 

Используя это наблюдение и определение функции $U(\alpha)$, сразу же получаем  первое неравенство в (\ref{equ.14}), так как
\[
U(\alpha)\le -1-\frac{1}{2\alpha}\biggl[-2\alpha-\frac{2\alpha^2}{1-2\alpha}\biggr]=\frac{\alpha}{1-2\alpha}
.
\]

Второе неравенство в (\ref{equ.14}) доказывается еще проще в силу того, что
$\log(1+x)\ge x-x^2/2$. \quad $\blacksquare$ 

\bigskip

В дальнейшем нам потребуется следующий технический  результат. 
\begin{lemma}\label{lemma.4} Пусть $\{p_k, \ k=1,\ldots,K-1\}$ и $\{q_k,\ k=1,\ldots\}$
-- неотрицательные последовательности, причем
\[
q_1=1,\ q_k\le \exp[-\rho (k-2)-1],\ k=1,\ldots, \quad \gamma>0.
\]
Обозначим 
\[
w_k=\frac{p_k}{P+Q}\mathbf{1}\{k<K\}+\frac{q_{k-K+1}}{P+Q}\mathbf{1}\{k\ge K\},
\]
где 
\[
P=\sum_{k=1}^K p_k , \quad  Q=\sum_{k=1}^\infty q_k\ge 1.
\]
Тогда 
\[
H(w)\stackrel{\rm def}{=}\sum_{k=1}^\infty w_k\log\frac{1}{w_k} \le \log\Bigl[K-1+{\rm e}^{R(\rho)}\Bigr],
\]
где
\[
R(\rho)=
\frac{2}{{\rm e}\rho}\mathbf{1}\{{\rm e}\rho<1\}+\biggl(1+\frac{1}{\rho{\rm e}}\biggr)
\exp\biggl(\frac{1-\rho{\rm e}}{1+\rho{\rm e}}\biggr)\mathbf{1}\{{\rm e}\rho\ge 1\}.
\]
\end{lemma}
\textit{Доказательство.} Используя выпуклость функции $\log(x)$, имеем
\begin{equation}\label{equ.15}
\begin{split}
H(w)=\frac{P}{P+Q}\sum_{k=1}^{K-1} \frac{p_k}{P}\log\frac{(P+Q)/P}{p_k/P}\\
+\frac{Q}{P+Q}\sum_{k=1}^\infty \frac{q_k}{Q}\log\frac{(P+Q)/Q}{q_k/Q}\\
\le \frac{P}{P+Q}\log\frac{P+Q}{P}+\frac{Q}{P+Q}\log\frac{P+Q}{Q}\\+\frac{P}{P+Q}\log(K-1)+\frac{1}{P+Q} \biggl[\sum_{k=1}^\infty q_k\log\frac{1}{q_k}+\log(Q)\sum_{i=1}^\infty q_k\biggr].
\end{split}
\end{equation}
Заметим, что функция $x\log(1/x)$ является монотонной при $x\in[0,{\rm e}^{-1}]$. Поэтому
\begin{equation}\label{equ.16}
\sum_{k=1}^\infty q_k\log\frac{1}{q_k}\le \sum_{k=2}^\infty
\exp[-\rho(k-2)-1][\rho(k-1)+1]\le \frac{2}{\rho {\rm e}}
\end{equation}
и также очевидно, что
\begin{equation}\label{equ.17}
\begin{split}
\log(Q)\sum_{k=1}^\infty q_k\le [\log(Q)]_+ \bigg\{1+\sum_{k=2}^\infty
\exp[-\rho(k-2)-1]\biggr\}\\
\le [\log(Q)]_+\biggl[\frac{1}{\rho {\rm e}}+1\biggr].
\end{split}
\end{equation}
Обозначим далее 
\[
x=\frac{Q}{P+Q}.
\]
Тогда из (\ref{equ.15}--\ref{equ.17}) находим
\begin{equation}\label{equ.18}
\begin{split}
H(w)\le \max_{x\in [0,1]}\biggl\{-x\log(x)-(1-x)\log(1-x)+(1-x)\log(K-1)\\
+
\frac{x}{{\rm e}\gamma}\max_{Q\ge 1} \frac{2+[\log(Q)]_+(1+\rho{\rm e})}{Q}\biggr\}.
\end{split}
\end{equation}
Нетрудно проверить, что при ${\rm e}\rho \le 1$
\begin{equation*}\label{en5}
\max_{Q\ge 1} \frac{2+[\log(Q)]_+(1+\rho{\rm e})}{Q}=2,
\end{equation*}
а при  ${\rm e}\rho > 1$
\begin{equation*}\label{en6}
\max_{Q\ge 1} \frac{2+[\log(Q)]_+(1+\rho{\rm e})}{Q}=(1+\rho{\rm e})
\exp\biggl(\frac{1-\rho{\rm e}}{1+\rho{\rm e}}\biggr).
\end{equation*}
Поэтому из (\ref{equ.18}) имеем
\begin{equation}\label{equ.19}
\begin{split}
H(w)\le \max_{x\in [0,1]}\biggl\{-x\log(x)-(1-x)\log(1-x)+(1-x)\log(K-1)\\
+
xR(\rho)\biggr\}.
\end{split}
\end{equation}
Легко проверить, что точка $x^*$, в которой достигается максимум в правой части этого неравенства, находится из уравнения
\begin{equation*}\label{en8}
\log{\frac{1-x^*}{x^*}}=\log(K-1)-R(\rho)
\end{equation*} 
и, следовательно, 
\[
x^*=\frac{1}{1+(K-1)\exp[-R(\rho)]}
\]
Поэтому из (\ref{equ.19})
\begin{equation*}
\begin{split}
H(w)\le\log(K-1)-\log(1-x^*)-x^*\biggl[\log\frac{x^*}{1-x^*}+\log(K-1)-R(\rho)\biggr]\\
=\log(K-1)-\log(1-x^*)=\log\Bigl[K-1+{\rm e}^{R(\rho)}\Bigr].
\end{split}
\end{equation*}
$\blacksquare$

\bigskip

\textit{Доказательство теоремы \ref{th.3}.} Положим
\begin{equation}\label{equ.20}
\begin{split}
M_\epsilon=\max \Bigl\{m: [\bar{r}(Y,\hat{\mu}^m)- r^{\mathcal{M}}(\mu)]\le 4\epsilon\sigma^2[m-\hat{m}(Y)]+4\sigma^2\Bigr\};
\end{split}
\end{equation}
здесь $\epsilon\in (0,1)$ -- некоторое положительное  число, которое будет выбрано позднее.

Смысл  введения этой величины состоит в том, что она позволяет разделить все несмещенные оценки рисков $\{\bar{r}(Y,\hat{\mu}^m), \ k=1,2,\ldots\}$ на два подмножества, ассоциированные с $\{m\le M_\epsilon\}$ и $\{m> M_\epsilon\}$. На первом подмножестве  индексов поведение несмещенных оценок $\bar{r}(Y,\hat{\mu}^m)$ носит в целом хаотичный характер, на втором  же подмножестве случайность наблюдается в существенно меньшей мере поскольку 
\begin{equation*}
\bar{r}(Y,\hat{\mu}^m)> \bar{r}^{\mathcal{M}}(\mu)+4\epsilon\sigma^2[m-\hat{m}(Y)]+4\sigma^2, \quad m>M_\epsilon.
\end{equation*}
Это свойство позволяет применить  лемму \ref{lemma.4}.  
Действительно, как  при доказательстве теоремы \ref{th.2} в \cite{LB}, логарифмируя ${w}^*_m(Y)$, находим
\begin{equation*}
\begin{split}
4\sigma^2\log \frac{1}{{w}^*_m(Y)}=\bar{r}(Y,\hat{\mu}^m)-\bar{r}^{\mathcal{M}}(Y)\\ +4\sigma^2
\log\biggr\{\sum_{s\in \mathcal{M}} \exp\biggl[-\frac{\bar{r}(Y,\hat{\mu}^s)- \bar{r}^{\mathcal{M}}(Y)}{4\sigma^2}\biggr]\biggr\}.
\end{split}
\end{equation*}
Далее в силу того, что $\hat{m}(Y)=\arg\min_m{\bar{r}(Y,\hat{\mu}^m})\le M_\epsilon$
\begin{equation*}
\log\biggr\{\sum_{s\in \mathcal{M}} \exp\biggl[-\frac{\bar{r}(Y,\hat{\mu}^s)- \bar{r}^{\mathcal{M}}(Y)}{4\sigma^2}\biggr]\biggr\}\ge 0.
\end{equation*}
Поэтому, применяя лемму \ref{lemma.4} и пользуясь выпуклостью функции $\log(\cdot)$, находим
\begin{equation}\label{ee1}
\begin{split}
\mathbf{E}\sum_{m\in \mathcal{M}} w_m^*(Y)[\bar{r}(Y,\hat{\mu}^m)- \bar{r}^{\mathcal{M}}(Y)]\le4\sigma^2 \mathbf{E} \sum_{m\in \mathcal{M}} w_m^*(Y)\log\frac{1}{w_m^*(Y)}\\ 
\le 4\sigma^2\log\Bigl[\mathbf{E}M_\epsilon-1+{\rm e}^{R(\epsilon)}\Bigr].
\end{split}
\end{equation}

Чтобы оценить сверху $\mathbf{E}M_\epsilon$, воспользуемся определением (\ref{equ.20}).  Тогда из  (\ref{equ.2}) получаем 
\begin{equation*}
\begin{split}
M_\epsilon = \max\biggl\{m:\ \sum_{i=m+1}^\infty\mu_i^2+ \sigma^2(1-4\epsilon)m -2\sigma \sum_{i=m+1}^\infty\mu_i \xi_i-
\sigma^2 \sum_{i=1}^{m}(\xi_i^2-1)  \le \\
\le \sum_{i=\hat{m}(Y)+1}^\infty\mu_i^2 +\sigma^2(1-4\epsilon)\hat{m}(Y) -2\sigma \sum_{i=\hat{m}(Y)+1}^\infty\mu_i \xi_i-
\sigma^2 \sum_{i=1}^{\hat{m}(Y)}(\xi_i^2-1) +4\sigma^2
\biggr\}.
\end{split}
\end{equation*}
	Далее, зафиксировав некоторое число $\gamma\in(0,1)$, отсюда находим
\begin{equation*}
\begin{split}
M_\epsilon =\max\biggl\{m:\ (1-\gamma)\sum_{i=m+1}^\infty\mu_i^2 +\sigma^2(1-4\epsilon-\gamma)m\\ 
+\gamma \sum_{i=m+1}^\infty\mu_i^2-2\sigma \sum_{i=m+1}^\infty\mu_i \xi_i
-\sigma^2 \sum_{i=1}^{m}(\xi_i^2-1) +\gamma \sigma^{2} m  \\
\le (1+\gamma)\sum_{i=\hat{m}(Y)+1}^\infty\mu_i^2 +\sigma^2(1+\gamma)\hat{m}(Y)\\-\gamma \sum_{i=\hat{m}(Y)+1}^\infty\mu_i^2-2\sigma \sum_{i=\hat{m}(Y)+1}^\infty\mu_i \xi_i
\\-\sigma^2 \sum_{i=1}^{\hat{m}(Y)}(\xi_i^2-1) -(\gamma+4\epsilon) \sigma^{2} \hat{m}(Y)
+4\sigma^2\biggr\}.
\end{split}
\end{equation*}
Следовательно
\begin{equation*}
\begin{split}
M_\epsilon\le M_\epsilon' =\max\biggl\{m:\ \sigma^2(1-4\epsilon-\gamma)m 
\\+\min_{s\ge 1}\biggl[\gamma \sum_{i=s+1}^\infty\mu_i^2-2\sigma \sum_{i=s+1}^\infty\mu_i \xi_i\biggr]
+\sigma^2 \min_{s\ge 1}\biggr[\sum_{i=1}^{s}(1-\xi_i^2) +\gamma  s\biggr]  \\
\le (1+\gamma)\biggr[\sum_{i=\hat{m}(Y)+1}^\infty\mu_i^2 +\sigma^2\hat{m}(Y)\biggr]\\+\max_{s\ge 1}\biggl[(\gamma+4\epsilon) \sum_{i=s+1}^\infty\mu_i^2-2\sigma \sum_{i=s+1}^\infty\mu_i \xi_i\biggr]\\
+\sigma^2 \max_{s\ge 1}\biggl[\sum_{i=1}^{s}(1-\xi_i^2) -\gamma s\biggr]
+4\sigma^2\biggr\}.
\end{split}
\end{equation*}

Для контроля максимумов и минимумов по $s$ в этом выражении используем леммы \ref{lemma.2} и \ref{lemma.3}. Тогда находим
\begin{equation*}
(1-4\epsilon-\gamma)\sigma^2\mathbf{E}M_\epsilon \le (1+\gamma)
\mathbf{E}\biggr[\sum_{i=\hat{m}(Y)+1}^\infty\mu_i^2 +\sigma^2\hat{m}(Y)\biggr]+4\sigma^2\biggl[2+\frac{2}{\gamma}\biggr].
\end{equation*} 
Следовательно, выбрав $\gamma=\epsilon$, получаем
\begin{equation}\label{equ.22}
\begin{split}
\sigma^2\mathbf{E}M_\epsilon \le  \frac{1+\epsilon}{1-5\epsilon}
\mathbf{E}\biggr[\sum_{i=\hat{m}(Y)+1}^\infty\mu_i^2
 +\sigma^2\hat{m}(Y)\biggr]
+\frac{10\sigma^2}{(1-5\epsilon)\epsilon}.
\end{split}
\end{equation}

Чтобы оценить  математическое ожидание в правой части этого неравенства,
воспользуемся тем, что для любого фиксированного $m\in \mathcal{M}$ выполнено неравенство
\begin{equation*}
-\sum_{i=1}^{\hat{m}(Y)}Y_i^2 +2\sigma^2 \hat{m}(Y)\le -\sum_{i=1}^{m}Y_i^2 +2\sigma^2 m
\end{equation*} 
или, что эквивалентно,
\begin{equation*}
\begin{split}
\sum_{i=\hat{m}(Y)+1}^\infty\mu_i^2 +\sigma^2 \hat{m}(Y) -2\sigma \sum_{i=\hat{m}(Y)+1}^\infty\mu_i \xi_i-
\sigma^2 \sum_{i=1}^{\hat{m}(Y)}(\xi_i^2-1) \\
\le\sum_{i=m+1}^\infty\mu_i^2+ \sigma^2 m -2\sigma \sum_{i=m+1}^\infty\mu_i \xi_i  -
\sigma^2 \sum_{i=1}^{m}(\xi_i^2-1). 
\end{split}
\end{equation*}
Выбрав произвольное число $\gamma\in (0,1)$, перепишем это неравенство в следующем виде:
\begin{equation*}
\begin{split}
(1-\gamma)\biggl[\sum_{i=\hat{m}(Y)+1}^\infty\mu_i^2 +\sigma^2 \hat{m}(Y) \biggr]
\le\sum_{i=m+1}^\infty\mu_i^2+ \sigma^2 m\\ -2\sigma \sum_{i=m+1}^\infty\mu_i \xi_i  -
\sigma^2 \sum_{i=1}^{m}(\xi_i^2-1)\\
+\max_{s\ge 1}\biggl[2\sigma \sum_{i=s+1}^\infty\mu_i \xi_i-\gamma \sum_{i=s+1}^\infty\mu_i^2 \biggr]\\
+\sigma^2\max_{s\ge 1}\biggl[ \sum_{i=1}^s\ (\xi_i^2-1)-\gamma s \biggr].
\end{split}
\end{equation*} 
Далее, воспользовавшись леммами \ref{lemma.2} и \ref{lemma.3}, находим
\begin{equation*}
(1-\gamma)\mathbf{E}\biggl[\sum_{i=\hat{m}(Y)+1}^\infty\mu_i^2 +\sigma^2 \hat{m}(Y) \biggr]
\le\sum_{i=m+1}^\infty\mu_i^2+ \sigma^2 m +\sigma^2\biggl[2+\frac{3}{\gamma}\biggr].
\end{equation*}
и минимизируя правую часть этого неравенства  по $m$, получаем
\begin{equation}\label{equ.23}
\mathbf{E}\biggl[\sum_{i=\hat{m}(Y)+1}^\infty\mu_i^2 +\sigma^2 \hat{m}(Y) \biggr]
\le r^{\mathcal{M}}(\mu)+ \frac{\gamma r^{\mathcal{M}}(\mu)}{1-\gamma} +\frac{\sigma^2}{(1-\gamma)}\biggl[2+\frac{3}{\gamma}\biggr].
\end{equation}
Поскольку это неравенство справедливо для любого  $\gamma\in (0,1)$, возьмем 
$
\gamma=\epsilon
$ 
и подставив эту величину в (\ref{equ.23}), получим
\begin{equation*}\label{et.10}
\begin{split}
\mathbf{E}\biggl[\sum_{i=\hat{m}(Y)+1}^\infty\mu_i^2
+\sigma^2 \hat{m}(Y) \biggr]
\le \frac{r^{\mathcal{M}}(\mu)}{1-\epsilon}+ \frac{4\sigma^2}{(1-\epsilon)\epsilon}.
\end{split}
\end{equation*} 
Подставляя это неравенство в (\ref{equ.22})
приходим к неравенству
\begin{equation*}
\begin{split}
\sigma^2\mathbf{E}M_\epsilon\le \frac{(1+\epsilon)r^{\mathcal{M}}(\mu)}{(1-5\epsilon)(1-\epsilon)}+
\frac{15}{(1-5\epsilon)\epsilon}\\
\le r^{\mathcal{M}}(\mu)+\frac{7\epsilon r^{\mathcal{M}}(\mu)}{(1-6\epsilon)}+
\frac{15}{(1-6\epsilon)\epsilon}.
\end{split}
\end{equation*}
Оъединяя это неравенство и (\ref{ee1}), получаем 
\begin{equation*}
\begin{split}
\mathbf{E}\sum_{m\in \mathcal{M}} w_m^*(Y)[\bar{r}(Y,\hat{\mu}^m)- \bar{r}^{\mathcal{M}}(Y)]
\le 4\sigma^2 \log\biggl[\frac{r^{\mathcal{M}}(\mu)}{\sigma^2}+\frac{7\epsilon r^{\mathcal{M}}(\mu)}{(1-6\epsilon)\sigma^2}\\+
\frac{15}{(1-6\epsilon)\epsilon}+\exp\biggl(\frac{2}{{\rm e}\epsilon}\biggr)\biggr].
\end{split}
\end{equation*} 
Для завершения доказательства теоремы остается минимизировать правую часть по $\epsilon\in[0,1/7]$.
Заметим, что  
\begin{equation*}
\begin{split}
\frac{7\epsilon r^{\mathcal{M}}(\mu)}{(1-6\epsilon)\sigma^2}+
\frac{15}{(1-6\epsilon)\epsilon}+\exp\biggl(\frac{2}{{\rm e}\epsilon}\biggr)\\ \le \frac{r^{\mathcal{M}}(\mu)}{\sigma^2}\biggl\{49\epsilon +\frac{\sigma^2}{r^{\mathcal{M}}(\mu)}\biggl[\frac{105}{\epsilon } +\exp\biggl(\frac{2}{{\rm e}\epsilon}\biggr)\biggr]\biggr\}.
\end{split}
\end{equation*}
Поэтому, выбирая
\[
\Psi(r)=\min_{\epsilon\in[0,1/7]}\biggl\{49\epsilon +r\biggl[\frac{105}{\epsilon }+\exp\biggl(\frac{2}{{\rm e}\epsilon}\biggr)\biggr]\biggr\}
\]
и, применяя лемму  \ref{lemma.1}, завершаем доказательство неравенства (\ref{equ.6}).

 Ясно, что $\Psi(r)$ ограничена при $r\in[0,1]$. 
Нетрудно проверить также, что при  
$r\rightarrow 0$
\begin{equation*}
\begin{split}
\epsilon^*(r)=\arg\min_{\epsilon\in[0,1/7]}\biggl\{49\epsilon +r\biggl[\frac{105}{\epsilon }+\exp\biggl(\frac{2}{{\rm e}\epsilon}\biggr)\biggr]\biggr\}\\
 \approx \frac{2}{\rm e}\log^{-1}\frac{49}{r}\biggl[1
+\log\biggl(\frac{2}{\rm e}\log\frac{49}{r}\biggr)\biggr]
\end{split}
\end{equation*}
и 
\[
\Psi(r)=(1+o(1))\frac{98}{\rm e}\log^{-1}\frac{49}{r}.
\]
$\blacksquare$

\bigskip

В заключение, автор хотел бы выразить  искреннюю признательность рецензенту за полезные и конструктивные замечания, способствовавшие улучшению работы.

\end{document}